\input amstex
\documentstyle {amsppt}
\UseAMSsymbols \vsize 18cm \widestnumber\key{ZZZZZ}

\catcode`\@=11
\def\displaylinesno #1{\displ@y\halign{
\hbox to\displaywidth{$\@lign\hfil\displaystyle##\hfil$}&
\llap{$##$}\crcr#1\crcr}}
\def\ldisplaylinesno #1{\displ@y\halign{
\hbox to\displaywidth{$\@lign\hfil\displaystyle##\hfil$}&
\kern-\displaywidth\rlap{$##$} \tabskip\displaywidth\crcr#1\crcr}}
\catcode`\@=12

\refstyle{A}

\let \ul=\underline
\let \ti=\widetilde

 \font\srm=cmr10 at 7.5pt

\font\main=cmsy10 at 10pt

\font\smain=cmsy10 at 7.5pt \font\ssmain=cmsy10 at 5.625pt

\font \fin=lasy8 at 15.4 pt

\def \o{\mathop{\hbox{\main O}}\nolimits}
\def \so{\mathop{\hbox{\smain O}}\nolimits}
\def \sso{\mathop{\hbox{\ssmain O}}\nolimits}

\def \ind{\mathop{\hbox {\rm ind}}\nolimits}

\def \A{\mathop{\hbox{\main A}}\nolimits}

\def \J{\mathop{\hbox{\main J}}\nolimits}

\def \GL{\mathop{\hbox{\rm GL}}\nolimits}
\def \End{\mathop{\hbox{\rm End}}\nolimits}

\def \Hom{\mathop{\hbox{\rm Hom}}\nolimits}

\def \GL{\mathop{\hbox{\rm GL}}\nolimits}
\def \SL{\mathop{\hbox{\rm SL}}\nolimits}

\def \Stab{\mathop{\hbox{\rm Stab}}\nolimits}

\topmatter
\title Param\`etres de Langlands et alg\`ebres d'entrelacement \endtitle

\rightheadtext{Param\`etres de Langlands et alg\`ebres
d'entrelacement}
\author Volker Heiermann \endauthor

\address Laboratoire de math\'ematiques, Universit\'e Blaise-Pascal, Aubi\`ere, France \endaddress

\email heiermann\@math.univ-bpclermont.fr \endemail

\thanks Ce travail a b\'enifici\'e d'une aide de l'Agence
Nationale de la Recherche portant la r\'eference
ANR-08-BLAN-0259-02. L'auteur remercie le rapporteur inconnu pour
son travail soigneux qui lui a permis de corriger quelques erreurs
de frappe dont au moins une de nature plus importante.
\endthanks

\abstract  Let $G$ be a classical $p$-adic group and $(\psi
,\epsilon )$ the Langlands parameter of an irreducible
supercuspidal representation of a Levi subgroup of $G$. Using data
from $(\psi ,\epsilon )$, we determine explicitly the intertwining
algebra of the representation which is induced from the orbit of
the supercuspidal representation associated to $(\psi ,\epsilon)$.
\rm

\null\noindent{R\'ESUM\'E: } Soit $G$ un groupe $p$-adique
classique et $(\psi ,\epsilon )$ le param\`etre de Langlands d'une
repr\'esentation irr\'eductible cuspidale d'un sous-groupe de Levi
de $G$. Utilisant des donn\'ees de $(\psi, \epsilon )$, nous
d\'eterminons explicitement l'alg\`ebre d'entrelacement de la
repr\'esentation induite par l'orbite de la repr\'esentation
cuspidale associ\'ee \`a $(\psi ,\epsilon )$.
\endabstract

\endtopmatter
\document

Fixons un corps local non archim\'edien $F$. On note
$\vert\cdot\vert_F$ sa valeur absolue normalis\'ee, $q$ la
cardinalit\'e de son corps r\'esiduel, $\varpi _F$ un
g\'en\'erateur de l'id\'eal maximal de son anneau de valuation et
$W_F$ son groupe de Weil.

Le symbole $G$ d\'esigne le groupe des points $F$-rationnels d'un
groupe classique connexe d\'efini sur $F$. Nous entendons ici par
groupe classique la composante neutre du groupe des automorphismes
d'un $F$-espace vectoriel laissant invariant une forme
bilin\'eaire sym\'etrique ou symplectique. C'est ou un groupe
symplectique ou un groupe orthogonal, \'eventuellement non
d\'eploy\'e.

Notons $\widehat{G}$ le groupe dual de Langlands associ\'e \`a
$G$. C'est un groupe r\'eductif complexe connexe qui est
orthogonal pair, si $G$ est orthogonal pair, symplectique si $G$
est orthogonal impair, et orthogonal impair si $G$ est
symplectique.

Soit $M$ un sous-groupe de Levi standard de $G$. Il s'identifie
\`a un groupe de la forme $\GL _{k_1}(F) \times\cdots\times
\GL_{k_h}(F)\times H$, o\`u $H$ d\'esigne un groupe classique du
m\^eme type que $G$. Notons $\GL_k(F)^1$ le sous-groupe de $\GL
_k(F)$ form\'e des \'el\'ements de d\'eterminant de valeur absolue
$1$. Posons $M^1=\GL _{k_1}(F)^1\times\cdots \times\GL_{k_h}(F)^1
\times H$. On appellera caract\`ere non ramifi\'e de $M$ tout
caract\`ere complexe de $M$ qui se factorise par $M^1$. Le groupe
des caract\`eres non ramifi\'es forme un tore alg\'ebrique
complexe, et l'alg\`ebre des fonctions r\'eguli\`eres de cette
vari\'et\'e alg\'ebrique affine s'identifie \`a l'alg\`ebre de
groupe $\Bbb C[M/M^1]$. On \'ecrira $B$ pour $\Bbb C[M/M^1]$ et,
pour $m$ dans $M$, $b_m$ pour la fonction r\'eguli\`ere qui
associe \`a un caract\`ere non ramifi\'e $\chi $ la valeur $\chi
(m)$.

Soit $(\sigma ,E)$ une repr\'esentation irr\'eductible cuspidale
de $M$. Notons $\o $ l'ensemble des classes d'isomorphie de
repr\'esentations de la forme $\sigma\otimes\chi $ avec $\chi $
caract\`ere non ramifi\'e de $M$.

Il a \'et\'e montr\'e par J. Bernstein \cite{BD} que la
cat\'egorie des repr\'esentations irr\'educ-tibles lisses $Rep(G)$
de $G$ se d\'ecompose en un produit direct de sous-cat\'egories
pleines $Rep(G)_{\so }$ index\'ees par les classes de conjugaison
d'orbites $\o $ de repr\'esenta-tions irr\'eductibles cuspidales
de sous-groupes de Levi $M$ de $G$.

Posons $E_B=E\otimes_{\Bbb C}B$, et notons $\sigma _B:M\rightarrow
E_B$ la repr\'esentation de $M$ d\'efinie par $\sigma
_B(m)(e\otimes b)=\sigma (m)e\otimes bb_m$. Soit $P=MU$ le
sous-groupe parabolique standard de $G$ contenant le sous-groupe
form\'e par des matrices triangulaires sup\'erieures, et notons
$i_P^G$ le foncteur d'induction parabolique normalis\'e qui
pr\'eserve l'unitarit\'e. Alors, un autre r\'esultat de J.
Bernstein \cite{Ru} dit que la cat\'egorie $Rep(G)_{\so }$ est
isomorphe \`a la cat\'egorie des modules \`a droite sur
l'alg\`ebre d'entrelacement $\End _G(i_P^GE_B)$.

Remarquons que, pour tout caract\`ere non ramifi\'e $\chi $ de
$M$, l'application de sp\'ecialisation $sp_{\chi }:B \rightarrow
\Bbb C$, $b\mapsto b(\chi )$, induit canoniquement un morphisme
$M$-\'equivariant $(\sigma _B,E_B)\rightarrow (\sigma\otimes\chi
,E_{\chi })$ et un morphisme $G$-\'equivariant $(i_P^G\sigma
_B,i_P^G$ $E_B)\rightarrow (i_P^G(\sigma\otimes\chi ),i_P^GE_{\chi
})$.

Notre but ici est de d\'ecrire l'alg\`ebre $\End _G(i_P^GE_B)$
explicitement en fonction du param\`etre de Langlands $(\psi
,\epsilon )$ de $\sigma $.

En effet, le param\`etre de Langlands $(\psi ,\epsilon )$ de
$\sigma $ est connu gr\^ace au travail de C. Moeglin \cite{M} sur
les r\'esultat de J. Arthur, lorsque $G$ est symplectique ou
orthogonal impair et $F$ de caract\'eristique $0$. Dans le cas de
la composante neutre d'un groupe orthogonal pair, il faut \^etre
un peu plus prudent et, par ailleurs, passer par le groupe
orthogonal tout entier qui n'est pas connexe, les r\'esultats de
C. Moeglin n'\'etant \'enonc\'es que dans ce cas. (Plus de
d\'etails sont donn\'es \`a la fin de la section 1.) La
restriction \`a la caract\'eristique $0$ devrait \^etre inutile,
mais pour le moment les r\'esultats utilis\'es relatifs aux
param\`etres de Langlands ne sont disponibles qu'en
caract\'eristique $0$.

D'autre part, on a calcul\'e dans \cite{H} l'alg\`ebre
d'entrelacement d'une certaine sous-repr\'esentation
$i_P^GE_{B_{\sso }}$ de $i_P^GE_B$ et montr\'e que celle-ci est
isomorphe \`a une alg\`ebre de Hecke avec param\`etres (ou
plut\^ot au produit semi-direct d'une telle alg\`ebre avec un
groupe fini). Rappelons \cite{Ro} que la cat\'egorie $Rep(G)_{\so
}$ est encore isomorphe \`a la cat\'egorie des modules \`a droite
sur l'alg\`ebre $\End _G(i_P^GE_{B_{\sso }})$.

Dans ce papier, on exprime d'abord le r\'esultat de \cite{H} sur
l'alg\`ebre d'entrelace-ment de $i_P^GE_{B_{\sso }}$ en termes du
param\`etre de Langlands de $\sigma $. Le r\'esultat principal
concernant l'alg\`ebre d'entrelacement de $i_P^GE_B$ se trouve
alors dans la section {\bf 5}. Pour la commodit\'e du lecteur, on
a rappel\'e au d\'ebut de cette sections toutes les notations et
d\'efinitions introduites ult\'erieurement.

L'auteur remercie C. Moeglin ainsi que C. Jantzen pour avoir
discut\'e avec lui sur leurs r\'esultats respectifs, ainsi que
A.-M. Aubert pour quelques corrections stylistiques.

\null{\bf 1.} Soit $F$ un corps local non archim\'edien de
caract\'eristique $0$. Soit $H$ un groupe symplectique, la
composante connexe d'un groupe orthogonal impair ou un groupe
orthogonal (non connexe) d\'efini sur $F$. Notons $\widehat{H}$
son groupe dual. Un param\`etre de Langlands temp\'er\'e pour $H$
est un couple $(\psi,\epsilon )$, o\`u $\psi :W_F\times SL_2(\Bbb
C)\rightarrow \widehat{H}$ est un homomorphisme qui v\'erifie les
propri\'et\'es suivantes: la restriction de $\psi $ \`a $W_F$ est
un homomorphisme continu dont l'image est born\'ee et form\'ee
d'\'el\'ements semi-simples. La restriction de $\psi $ \`a
$\SL_2(\Bbb C)$ est un homomorphisme de groupes alg\'ebriques. La
deuxi\`eme composante, $\epsilon $, est un certain caract\`ere du
centralisateur de l'image de $\psi $ qui doit par ailleurs avoir
la "bonne" restriction au centre de $\widehat{H}$.

Le param\`etre est dit \it discret, \rm si l'image de $\psi $
n'est contenu dans aucun sous-groupe de Levi propre de $\widehat
{H}$.

\null {\bf 1.1} Il r\'esulte des travaux de J. Arthur - via la
correspondance de Langlands locale pour les groupes lin\'eaires
g\'en\'eraux $p$-adiques - qu'\`a tout param\`etre discret
$(\psi,\epsilon )$ pour $H$ correspond une unique repr\'esentation
de carr\'e int\'egrable $\tau (\psi,\epsilon )$ de $H$, et
vice-versa, comme cela a \'et\'e conjectur\'e par Langlands (avec
des raffinements ult\'erieurs de P. Deligne et G. Lusztig).

\null {\bf 1.2} C. Moeglin \cite{M} a su expliciter les
param\`etres qui correspondent \`a des repr\'e-sentations
cuspidales. Explicitons cela: notons $\iota
:\widehat{H}\rightarrow\GL _{k_{\widehat{H}}}(\Bbb C)$ la
repr\'esentation naturelle de $\widehat{H}$ et $Jord(\psi )$ le
multi-ensemble des composantes irr\'eductibles de $\iota\circ\psi
$. Ses \'el\'ements sont des repr\'esentations de la forme
$\rho\otimes sp_a: W_F\times SL_2(\Bbb C)\rightarrow \GL
_{ad_{\rho }}(\Bbb C)$, o\`u $\rho :W_F\rightarrow \GL _{d_{\rho
}}(\Bbb C)$ est une repr\'esentation irr\'eductible autoduale, $a$
un entier $\geq 1$, et $sp_a$ la repr\'esentation irr\'eductible
de $\SL _2(F)$ de degr\'e $a$. On a donc
$$\iota\circ\psi =\bigoplus _{(\rho ,a)\in Jord(\psi )}\rho\otimes
sp_a.$$ C. Moeglin a montr\'e \cite{M, 1.5} que $\tau (\psi
,\epsilon )$ est cuspidale, si et seulement si l'ensemble
$Jord(\psi )$ est sans trou (i.e. $(\rho ,a)\in Jord(\psi )$,
$a\geq 3$, implique $(\rho, a-2)\in Jord(\psi )$) et si le
caract\`ere $\epsilon $ est altern\'e avec la bonne restriction au
centre de $\widehat{H}$. (On ne va pas expliquer ici, ce que cela
signifie, puisque nous n'en aurons pas besoin dans la suite.
Signalons seulement que l'existence d'un caract\`ere $\epsilon $
n'est pas automatiquement v\'erifi\'ee, si $Jord(\psi )$ est sans
trou.)

Si $(\rho ,a)\in Jord(\psi )$, alors $a$ est impair si $\rho $ et
$\widehat {H}$ sont tous les deux orthogonaux ou tous les deux
symplectiques. Dans le cas contraire, $a$ est pair.

Si $\rho $ est une repr\'esentation irr\'eductible autoduale de
$W_F$, notons, si cet entier existe, $a_{\rho, \psi }$ le plus
grand entier $a$ tel que $(\rho ,a)\in Jord(\psi )$. Sinon, posons
$a_{\rho, \psi }=-1$ si $\rho $ et $\widehat{H}$ sont tous les
deux orthogonaux ou symplectiques, et $a_{\rho,\psi }=0$ sinon.

Notons encore $\rho $ la repr\'esentation irr\'eductible cuspidale
de $\GL _{d_{\rho }}(F)$ qui correspond \`a $\rho $ par la
correspondance locale de Langlands. Elle est autoduale. Si $\GL
_{d_{\rho }}(F)\times H$ est un sous-groupe de Levi de $G$, la
repr\'esentation de $G$ d\'eduite de $\rho \vert\cdot\vert
_F^x\otimes\tau (\psi )$ par induction parabolique normalis\'ee
est r\'eductible pour un seul r\'eel $x\geq 0$. C. Moeglin
\cite{M} a montr\'e que cet entier vaut $(a_{\rho ,\psi }+1)/2$.
(Lorsque $a_{\rho ,\psi }$ vaut $-1$ ou $0$, la condition donn\'ee
dans \cite{M} porte sur le transfert d'une certaine distribution
stable. L'\'equivalence avec la condition ci-dessus est une
cons\'equence des r\'esultats annonc\'es par J. Arthur.) Lorsque
$H$ est trivial, on consid\`ere $\widehat{H}$ suivant la nature de
$G$ comme orthogonal ou symplectique. \'Evidemment, lorsque par
exemple $H=1$ et $\rho $ est un caract\`ere autodual, ces
r\'esultats sont connus depuis longtemps.

\null{\bf 1.3 Proposition:} \it Soit $\rho $ une repr\'esentation
irr\'eductible autoduale de $W_F$. Soit $\chi _-$ un caract\`ere
non ramifi\'e de $W_F$ tel que la repr\'esentation $\rho
_-=\rho\chi _-$ soit autoduale et non isomorphe \`a $\rho $.
Notons $t_{\rho }$ l'ordre du groupe des caract\`eres non
ramifi\'es $\chi $ de $W_F$ v\'erifiant $\rho
\simeq\rho\otimes\chi $.

Les entiers $t_{\rho }a_{\rho ,\psi }$ et $t_{\rho _-}a_{\rho
_-,\psi }$ ont la m\^eme parit\'e.

\null Preuve: \rm Ces entiers sont impairs, si et seulement si
$t_{\rho }$ et $a_{\rho ,\psi }$ (resp. $t_{\rho _-}$ et $a_{\rho
_-,\psi }$) sont tous les deux impairs. Remarquons que $t_{\rho }$
et $t_{\rho _-}$ sont \'egaux. Comme $\rho $ et $\rho \chi _-$
sont autoduales, il est imm\'ediat que $\rho\chi _-^2\simeq\rho $.
Supposons $t_{\rho }$ impair. Quitte \`a multiplier $\chi _-$ par
un caract\`ere $\chi $ qui v\'erifie $\rho\simeq\rho\otimes\chi $
(ce qui ne change pas la classe d'isomorphie de $\rho _-$), on
peut supposer dans ce cas $\chi _-^2=1$.

Mais, alors $\rho $ est orthogonal (resp. symplectique), si et
seulement si $\rho _-$ l'est. Comme l'entier $a_{\rho ,\psi }$
(resp. $a_{\rho _-,\psi }$) est impair si et seulement si
$\widehat {H}$ et $\rho $ (resp. $\rho _-$) sont tous les deux
orthogonaux ou symplectiques, ceci prouve la proposition.
\hfill{\fin 2}

\null{\bf 1.4} Le symbole $H$ d\'esignera maintenant la composante
connexe d'un groupe orthogonal pair $H'$ d\'efini sur $F$. Le
r\'esultat suivant est la proposition {\bf 4.3} de \cite{BJ}
(compl\'et\'e au cas non d\'eploy\'e par \cite{J}):

\null{\bf Proposition:} \it Soit $\tau $ une repr\'esentation
irr\'eductible cuspidale de $H$, et soit $\tau '$ une composante
irr\'eductible de la repr\'esentation induite $ind_H^{H'}\tau $.
Soit $\rho $ une repr\'esentation irr\'eductible cuspidale de $\GL
_k(F)$. Notons $G'$ (resp. $G$) le groupe orthogonal (resp. sa
composante connexe) dont $\GL _k(F)\times H'$ (resp.
$\GL_k(F)\times H$) s'identifie \`a un sous-groupe de Levi
standard maximal. Notons $P'$ (resp. $P$) le sous-groupe
parabolique standard maximal de $G'$ (resp. $G$) de Levi $M'$
(resp. $M$). Notons $c$ un repr\'esentant de l'\'el\'ement non
trivial de $G'/G$.

(i) Supposons $G$ d\'eploy\'e, $k$ impair et ou bien
$\tau\not\simeq c\tau$ ou bien $H=1$ et $k\ne 1$. Alors
$i_P^G(\rho\vert\cdot\vert^x\otimes\tau )$ est irr\'eductible pour
tout nombre r\'eel $x$.

(ii) Supposons les conditions de (i) non v\'erifi\'ees. Fixons un
nombre r\'eel $x$. Alors, pour que
$i_P^G(\rho\vert\cdot\vert^x\otimes\tau )$ soit r\'eductible, il
faut et il suffit que
$i_{P'}^{G'}(\rho\vert\cdot\vert^x\otimes\tau ')$ soit
r\'eductible. \rm

\null Cette proposition nous conduit \`a introduire la
terminologie suivante: appelons \it param\`etre de Langlands \rm
d'une repr\'esentation irr\'eductible cuspidale $\tau $ de $H$, le
param\`etre de Langlands $(\psi ,\epsilon )$ d'une composante
irr\'eductible de la repr\'esentation induite $ind_H^{H'}\tau $.
Ceci est \'evidemment un abus de terminologie qui nous semble
toutefois justifi\'e par le corollaire suivant qui est une
cons\'equence imm\'ediate de la proposition ci-dessus et des
r\'esultats de C. Moeglin dans le cas d'un groupe orthogonal pair
(non connexe) d\'ecrits dans la premi\`ere partie de cette
section:

\null{\bf Corollaire:} \it Soit $(\psi ,\epsilon )$ un param\`etre
de Langlands d'une repr\'esentation irr\'educ-tible cuspidale
$\tau $ de $H$. Soit $\rho $ une repr\'esentation irr\'eductible
cuspidale autoduale de $\GL _k(F)$. Notons $G$ un groupe
orthogonal connexe dont $\GL_k(F)\times H$ s'identifie \`a un
sous-groupe de Levi standard maximal.

Alors, la repr\'esentation induite de
$\rho\vert\cdot\vert_F^x\otimes\tau $ par induction parabolique
normalis\'ee est toujours irr\'eductible si $G$ est d\'eploy\'e,
$k$ impair et ou bien $H=1$ et $k\ne 1$ ou bien
$c\tau\not\simeq\tau $. Sinon, elle est r\'eductible pour un seul
nombre r\'eel $x\geq 0$, et celui-ci vaut $(a_{\rho ,\psi }+1)/2$.
\rm

\null{\bf 2.} Soit maintenant $G$ un groupe symplectique ou la
composante connexe d'un groupe orthogonal. On appellera \it
param\`etre de Langlands d'une repr\'esentation cuspidale d'un
sous-groupe de Levi de $G$ \rm un couple $(\psi ,\epsilon )$,
form\'e d'un homomorphisme $\psi :W_F\times\SL_2(\Bbb
C)\rightarrow\widehat{G}$ qui v\'erifie par ailleurs la
propri\'et\'e suivante: si $M$ est un sous-groupe de Levi de $G$
qui est minimal pour la propri\'et\'e que l'image de $\psi $ est
contenue dans le groupe dual $\widehat{M}$, alors $(\psi ,\epsilon
)$ est le param\`etre d'une repr\'esentation cuspidale $(\sigma
,E)$ de $M$. On dira alors que $M$ est un sous-groupe de Levi
associ\'e \`a $(\psi ,\epsilon )$.

Deux tels param\`etres seront dits \'equivalents, s'ils sont
conjugu\'es par un \'el\'ement de $\widehat{G}$. Remarquons que la
classe d'isomorphie de la repr\'esentation $i_P^GE_B$ ne change
pas, si on passe \`a un param\`etre \'equivalent. (Dans le cas de
la composante connexe d'un groupe orthogonal pair, cette notion se
g\'en\'eralise de fa\c con \'evidente.)

\null{\bf 2.1} Soit $M$ un sous-groupe de Levi standard de $G$ et
$(\psi ,\epsilon )$ le param\`etre de Langlands d'une
repr\'esentation irr\'eductible cuspidale de $M$. Quitte \`a
conjuguer $(\psi ,\epsilon )$ par un \'el\'ement de $\widehat{G}$,
on peut toujours supposer que $\psi $ soit de la forme $$\ti{\rho
}_{1,1}\otimes\cdots\otimes\ti{\rho }_{1,d_1}\otimes\ti{\rho
}_{2,1}\otimes\cdots\otimes\ti{\rho
}_{2,d_2}\otimes\cdots\otimes\ti{\rho }_{h,d_h}\otimes\psi _H,$$
o\`u les $\ti{\rho }_{i,j}:W_F\rightarrow\GL _{k_i}(\Bbb C)$ sont
des repr\'esentations de $W_F$ qui sont, pour $i$ fix\'e, la
tordue d'une m\^eme repr\'esentation irr\'eductible cuspidale
unitaire $\rho _i$ par un caract\`ere non ramifi\'e. On peut
choisir (et on choisira) $\rho _i$ autoduale si la
contragr\'ediente de $\ti{\rho }_{i,j}$ est isomorphe au produit
de $\ti{\rho }_{i,j}$ par un caract\`ere non ramifi\'e, et on
supposera que $\rho _i$ n'est pas la tordue d'une $\rho _j$, $j\ne
i$, par un caract\`ere non ramifi\'e. Le couple $(\psi _H,\epsilon
)$ est le param\`etre de Langlands d'une repr\'esentation
irr\'eductible cuspidale $(\tau ,E_{\tau })$ d'un groupe $H$ du
m\^eme type que $G$, mais de rang plus petit, \'eventuellement
trivial. Notons $\rho _{i-}$ la repr\'esentation de $W_F$,
d\'etermin\'ee \`a isomorphisme pr\`es, qui est le produit de
$\rho _i$ par un caract\`ere non ramifi\'e et qui est autoduale et
non isomorphe \`a $\rho _i$. On peut par ailleurs supposer (et on
supposera) $a_{\rho _i,\psi _H}\geq a_{\rho _{i-},\psi _H}$. (Ce
choix des $\rho _i$ est conforme au choix du point de base
effectu\'e dans \cite{H}.) Le sous-groupe de Levi $M$ s'identifie
alors \`a un produit $$\GL
_{k_1}(F)\times\cdots\times\GL_{k_1}(F)\times\GL_{k_2}(F)\times
\cdots \times\GL_{k_h}\times\cdots\times\GL_{k_h}(F)\times H,$$
chaque facteur $\GL _{k_i}$ \'etant r\'ep\'et\'e $d_i$ fois.

\null {\bf 2.2} Plus pr\'ecis\'ement, dans la r\'ealisation
usuelle de $G$ comme sous-groupe de $\GL_{k_{\widehat{G}}}(F)$,
$M$ est l'ensemble des \'el\'ements du groupe $\GL _{k_1}(F)\times
\cdots \times \GL_{k_h}(F)\times H\times\GL _{k_h}(F)\times \cdots
\times \GL_{k_1}(F)$, plong\'e diagonalement dans $\GL
_{k_{\widehat{G}}} (F)$, qui sont de la forme
$(m_{1,1},m_{1,2},\dots ,m_{1,d_1}, m_{2,1},\dots
,m_{h,d_h},m_H,^tm_{h,d_h}^{-1},\dots ,^tm_{1,1}^{-1})$.

Pour $i=1,\dots ,h$ et $j=1,\dots ,d_i-1$, on notera $r_{i,j}$
l'\'el\'ement du groupe de Weyl de $G$ dont l'action sur $M$
permute les coefficients $m_{i,j}$ et $m_{i,j+1}$ (ainsi que
$^tm_{i,j}^{-1}$ et $^tm_{i,j+1}^{-1}$) d'un \'el\'ement de $M$.
Sauf si $G$ est orthogonal pair, $k_i$ impair et ou bien $H=1$ ou
bien $\tau $ n'est pas stable par l'automorphisme ext\'erieur, on
d\'esignera par ailleurs, si $\rho _i$ est autoduale, par
$r_{i,d_i}$ l'\'el\'ement du groupe de Weyl de $G$ dont l'action
sur $M$ permute les coefficient $m_{i,d_i}$ et $^tm_{i,d_i}^{-1}$
de $M$.

On \'ecrira $W_{\psi }$ pour le sous-groupe du groupe de Weyl de
$G$ engendr\'e par les $r_{i,j}$, et $W_{\psi ,i}$ pour le
sous-groupe engendr\'e par les $r_{i,j}$ avec $i$ fix\'e. Le
groupe $W_{\psi }$ est le produit direct des $W_{\psi ,i}$. On
identifie les \'el\'ements de $W_{\psi }$ \`a des \'el\'ements de
$G$ \`a l'aide d'un choix de repr\'esentants dans un certain
sous-groupe compact de $G$.

\null {\bf 2.3} Notons par abus de notations encore $(\rho
_i,E_{\rho _i})$ la repr\'esentation irr\'eductible cuspidale de
$\GL _{k_i}(F)$ qui correspond \`a $\rho _i$. Posons $$\sigma=\rho
_1\otimes\cdots\otimes \rho_1\otimes\rho
_2\otimes\cdots\otimes\rho _h\otimes\tau ,$$ chaque facteur $\rho
_i$ \'etant r\'ep\'et\'e $d_i$ fois. Notons $E$ l'espace de cette
repr\'esentation. Ce n'est en g\'en\'eral pas la repr\'esentation
de $M$ qui correspond au param\`etre de Langlands $(\psi ,\epsilon
)$, mais le produit de celle-ci par un certain caract\`ere non
ramifi\'e. Ce choix ne change toutefois pas la classe d'isomorphie
de la repr\'esentation $i_P^GE_B$.

\null {\bf 2.4} Rappelons que $\GL_{k_i}(F)^1$ d\'esigne le
sous-groupe de $\GL_{k_i}(F)$ form\'e des \'el\'ements de
d\'eterminant de valeur absolue $1$. Le groupe quotient
$\GL_k(F)/\GL_k(F)^1$ est cyclique, engendr\'e par l'image de la
matrice diagonale $h_{k_i}:=diag(\varpi _F, 1, 1,\dots ,$ $1)$.
Notons $Stab(\rho _i)$ le groupe des caract\`eres non ramifi\'es
$\chi $ de $\GL_{k_i}(F)$ (i.e. de restriction triviale \`a
$\GL_{k_i}(F)^1$) tels que $\rho _i\otimes\chi $ soit isomorphe
\`a $\rho _i$, et $t_i$ l'ordre de $Stab(\rho _i)$. Il r\'esulte
de la correspondance de Langlands que ce nombre $t_i$ est \'egal
au nombre $t_{\rho _i}$ d\'efini en {\bf 1.3} relatif \`a la
repr\'esentation galoisienne $\rho _i$.

Fixons une composante irr\'eductible $\rho _i^1$ de $(\rho
_i)_{\vert \hbox{\srm GL}_k(F)^1}$. Notons $E_{\rho _i}^1$ le
sous-espace de $E_{\rho _i}$ correspondant \`a $\rho _i^1$.

\null{\bf Proposition:} \it (cf. \cite{H,1.16}) On a $$E_{\rho
_i}=\oplus _{j=1}^{t_i}\ \rho _i(h_{k_i}^j)E_{\rho _i}^1.$$

\null \rm Posons $E_{\rho _i}^j=\rho _i(h_{k_i}^j)E_{\rho }^1$,
d\'efinissons $\J =\displaystyle{\times_{i=1}^h}\{1,\dots
,t_i\}^{d_i}$ et, pour $\underline{j}=(j_{1,1},\dots $
$,j_{h,d_h})$ dans $\J $, posons $E^{\underline{j}}=E_{\rho
_1}^{j_{1,1}}\otimes\cdots\otimes E_{\rho _h}^{j_{h,d_h}}\otimes
E_{\tau }$. Ce sont des sous-espaces irr\'eductibles de $E_{\rho
_i}$ et $E$ pour les actions de $\GL _{k_i}^1(F)$ et $M^1$
respectivement. Notons les repr\'esentations correspondantes
respectivement $\rho _i^j$ et $\sigma ^{\underline{j}}$. Elles
sont deux \`a deux in\'equivalentes (cf. \cite{H, 1.16}), et on a
$E=\bigoplus_{\underline{j}\in\Cal J}E^{\underline{j}}$.

Rappelons que toute sym\'etrie simple $r_{i,j}$ de $W_{\psi }$
v\'erifie $r_{i,j}\sigma\simeq\sigma $. On a donc $(r_{i,j}\sigma
^{\underline{j}},r_{i,j}E^{\underline{j}})=(\sigma
^{r_{i,j}(\underline{j})},E^{r_{i,j}(\underline{j}}))$ pour un
certain $r_{i,j}(\underline{j})$ qui se calcule de la mani\`ere
suivante: si $j=1,\dots ,d_i-1$, alors $r_{i,j}(\underline{j})$ se
d\'eduit de $\underline{j}$, en \'echangeant $j_{i,j}$ et
$j_{i,j+1}$. Si $j=d_i$, alors $\rho _i$ est autoduale, et il
existe $j_{i,d_i}'$ tel que la contragr\'ediente $(E_{\rho
_i}^{j_{i,d_i}})^{\vee }$ de $E_{\rho _i}^{j_{i,d_i}}$ soit
isomorphe \`a $E_{\rho _i}^{j_{i,d_i}'}$. Alors, $r_{i,d_i}
(\underline{j})$ se d\'eduit de $\underline{j}$, en rempla\c cant
$j_{i,d_i}$ par $j_{i,d_i}'$.

\null{\bf 3.} Dans cette section pr\'eliminaire, on part de la
situation donn\'ee dans {\bf 2.1}. On se fixe un entier $i$
compris entre $1$ et $h$, et on \'ecrira $\rho =\rho _i$,
$\rho_-=\rho_{i-}$, $k=k_i$, $r_j=r_{i,j}$, $d=d_i$, $t=t_i$ etc.

\null{\bf 3.1} On distinguera dans la suite trois cas. Supposons
d'abord $G$ symplectique ou orthogonal impair. Alors ces trois cas
sont

\null (I) la repr\'esentation $\rho $ n'est pas autoduale;

(II) la repr\'esentation $\rho $ est autoduale, la
repr\'esentation galoisienne $\rho $ ne figure pas dans $Jord(\psi
_H)$, et $\rho $ et $\widehat{G}$ sont tous les deux ou
orthogonaux ou symplectiques;

(III) la repr\'esentation $\rho $ est autoduale et la
repr\'esentation galoisienne $\rho $ figure ou dans $Jord(\psi
_H)$ ou $\rho $ et $\widehat{G}$ ne sont pas de la m\^eme nature
(i.e. l'un est symplectique et l'autre orthogonal).

\null \noindent{Si $G$ est orthogonal pair et d\'eploy\'e, alors
on ajoute \`a (I) le cas $k$ impair avec ou bien $H=1$ et $k\ne 1$
ou bien $\tau(\psi _H)$ non invariant par l'automorphisme
ext\'erieur, et on l'enl\`eve des autres cas. Par ailleurs, on
introduit le cas suivant:

\null (IIb) $G$ est orthogonal pair et d\'eploy\'e, $H=1$, $k=1$
et la repr\'esentation $\rho $ est autoduale.}

\null Posons $d'=d-1$ dans le cas (I), $d'=d$ dans les deux autres
cas. De plus, $s_j:=r_j$ pour $j=1,\dots ,d-1$,
$s_d=r_dr_{d-1}r_d^{-1}$ dans les cas (II) et (IIb), et $s_d=r_d$
dans le cas (III). (Si l'entier $i$ n'est plus fix\'e, on \'ecrira
$s_{i,j}$.) Si $d=1$, on laisse $s_1$ ind\'efini dans les cas (I),
(II) et (IIb).

\null{\bf 3.2} Pour $j=1,\dots ,d'$ et $\chi $ un caract\`ere non
ramifi\'e de $M$, notons $J_{s_jP\vert P}(\sigma\otimes\chi )$
l'op\'erateur d'entrelacement standard tel que d\'efini dans
\cite{W, IV}. Cet op\'erateur est d\'efini pour $\chi $ en dehors
d'un nombre fini de hyperplans. Il entrelace alors les
repr\'esentations $i_P^G(\sigma \otimes\chi )$ et
$i_{s_jP}^G(\sigma\otimes \chi )$. Il existe un \'el\'ement
$\ti{J}_{s_j}$ de $\Hom _G(i_P^GE_B,$ $i_{s_jP}^G E_B)$ et un
\'el\'ement $p_{s_j}\in B$ tel que $sp_{\chi
}\ti{J}_{s_j}=p_{s_j}(\chi )J_{s_jP\vert P}(\sigma\otimes\chi
)sp_{\chi }$. (L'homo-morphisme de sp\'ecialisation $sp_{\chi }$ a
\'et\'e d\'efini dans l'introduction.)

Notons $\lambda $ l'op\'eration par translations \`a gauche de $G$
sur $i_P^GE_B$ (le groupe $G$ agissant par translations \`a
droite) et $\tau _{s_j}$ l'automorphisme de $B$ qui envoie $b_m$
sur $b_{s_j^{-1}ms_j}$. Suivant \cite{H,2.4}, on d\'efinit de la
mani\`ere suivante, pour $j=1,\dots ,d-1$, un isomorphisme $\rho
_{s_j}:i_P^Gs_jE\rightarrow i_P^GE$: notons $h_k^{(i,j)}$ (ou
$h_k^{(j)}$ si $i$ est fix\'e) l'\'el\'ement de $M$ dont toutes
les entr\'ees sont \'egales \`a $1$, sauf celle en position
$(i,j)$ qui vaut $h_k$. Si $j$ est compris entre $1$ et $d-1$, on
pose
$$\rho _{s_j}=[(\chi(h_k^{(j)}{h_k^{(j+1)}}^{-1})-1)\lambda (s_j)
J_{s_j^{-1}P\vert P} (\sigma\otimes\chi )]_{\vert \chi
=1}]^{-1}.$$ Dans les cas (II) et (IIb), on pose, si $d>1$, $$\rho
_{s_d}=[(\chi(h_k^{(d)})-1) (\chi(h_k^{(d-1)}{h_k^{(d)}})-1)
(\chi(h_k^{(d-1)})-1)\lambda (s_d) J_{s_d^{-1}P\vert P}
(\sigma\otimes\chi))_{\vert \chi =1}]^{-1},$$ et dans le cas (III)
$$\rho _{s_d}=[(\chi(h_k^{(d)})-1)\lambda (s_d) J_{s_d^{-1}P\vert P}
(\sigma\otimes\chi )]_{\vert \chi =1}]^{-1}.$$ D\'esignons par
$K(B)$ le corps des fractions de $B$. L'op\'erateur
$A_{s_j}=\rho_{s_j}\tau_{s_j}\lambda (s_j)p_{s_j}^{-1}$
$\ti{J}_{s_j}$ est un \'el\'ement de $\Hom _G(i_P^GE_B,i_P^G
E_{K(B)})$ (cf. \cite{H,3.1}).

Dans le cas (II), on se fixe par ailleurs un isomorphisme $\rho
_{r_d}:r_d\sigma\rightarrow\sigma $ prolong\'e par fonctorialit\'e
en un isomorphisme $i_P^Gr_d\sigma\rightarrow i_P^G\sigma $ qui
sera toujours not\'e $\rho _{r_d}$, et on \'ecrit $A_{r_d}=\rho
_{r_d}\tau_{r_d}\lambda(r_d)p_{r_d}^{-1}\ti{J}_{r_d}$. C'est
\'egalement un \'el\'ement de $\Hom _G(i_P^GE_B,i_P^G$ $E_{K(B)})$
(cf. \cite{H,3.1}).

\null {\bf 3.3} D\'esignons par $ind_{M^1}^M$ le foncteur de
l'induction compacte.  Rappelons que $E_B=ind_{M^1}^M E_{\vert
M^1}$. On va d'abord s'int\'eresser au sous-espace
$i_P^G(ind_{M^1}^ME^{\ul{j}})$ de $i_P^GE_B$. L'op\'erateur
$A_{s_j}$ ne laisse en g\'en\'eral pas stable ce sous-espace. Pour r\'em\'edier \`a
cela, il faut le multiplier avec un certain \'el\'ement $b\in
B^{\times }$. La proposition ci-dessous est une cons\'equence
imm\'ediate de \cite {H, 4.5} et de sa preuve.

\null{\bf Proposition:} \it Soit $\underline{j}=(j_{1,1},\dots
,j_{h,d_h})\in\J $. Pour $1\leq j\leq d-1$, posons
$h_{\underline{j},s_j}=(h_k^{(i,j)}{h_k^{(i,j+1)}}^{-1}
)^{j_{i,j}-j_{i,j+1}}$. D\'efinissons
$h_{\ul{j},s_d}=(h_k^{(i,d-1)})^{j_{i,d-1}-j_{i,d}'}$
$(h_k^{(i,d)})^{j_ {i,d}-j_{i,d-1}'}$ dans les cas (II) et (IIb)
et posons $h_{\underline{j},s_d}=(h^{(i,d)})^{j_{i,d} -j_{i,d}'}$
dans le cas (III). Ce dernier \'el\'ement sera not\'e
$h_{\ul{j},r_d}$ dans le cas (II).

L'op\'erateur $b_{h_{\underline{j},s_j}}A_{s_j}$ laisse invariant
l'espace $i_P^G(ind_{M^1}^M E^{\underline{j}})$. Dans le cas (II),
il en est de m\^eme de l'op\'erateur $b_{h_{\ul{j},r_d}}A_{r_d}$.
\rm

\null{\bf 3.4} \'Ecrivons $J^{\underline{j}}_{s_j}$ pour
l'\'el\'ement de $\End _G(i_P^G(ind_{M^1}^M E^{\underline{j}}))$
qui est \'egal \`a la restriction de
$b_{h_{\underline{j},s_j}}A_{s_j}$, ainsi que, dans les cas (II),
$J_{r_d}^{\ul{j}}$ pour la restriction de
$b_{h_{\ul{j},r_d}}A_{r_d}$.

Pour $j=1,\dots ,d-1$, posons
$X_j=b_{h_k^{(j)}{h_k^{(j+1)}}^{-1}}^t$,
$X_d=b_{h_k^{(d-1)}h_k^{(d)}}^t$ dans les cas (II) et (IIb) si
$d>1$, et $X_{d}=b_{h_k^{(d)}}^t$ dans le cas (III). D\'efinissons
par ailleurs $a=a_{\rho ,\psi }$ et $a_-=a_{\rho _-,\psi }$ (cf.
{\bf 1.2}).

D'apr\`es \cite{H, 5.2}, il existe, pour $j=1,\dots ,d-1$ ou bien
encore $j=d$ dans les cas (II) et (IIb), un scalaire $c$ tel que
$sp_1(X_j-1)J_{s_j}=c(1-q^{-t})sp_1$. Dans le cas (III), fixons un
caract\`ere non ramifi\'e $\chi _-$ tel que $\rho\otimes\chi _-$
soit isomorphe \`a $\rho _-$. Il existe alors des scalaires $c_-$
et $c_+$ qui ne diff\`erent que par un facteur $\pm 1$ tels que
$$sp_1(X_d-1)J_{s_d}=c_+{(1-q^{-t{a+1\over 2}}) (1+q^{-t{a_-+1\over
2}})\over 2}sp_1$$
$$sp_{\chi _-}(X_d+1)J_{s_d}=c_-{(1+q^{-t{a+1\over 2}}) (1-q^{-t{a_-+1\over
2}})\over 2}sp_{\chi _-}.\leqno{\hbox{\rm et}}$$ Posons, pour $j=1,\dots ,d-1$,
$$T_{s_j}^{\underline{j}}= q^t cJ_{s_j}^{\underline{j}}-
(q^t-1){X_j\over 1-X_j},$$ dans les cas (II) et (IIb)
$$T_{s_d}^{\underline{j}}=q^tc J_{s_d}^{\underline{j}}- (q^t-1)
{X_d\over 1-X_d},$$ et dans le cas (III), lorsque
$c_+a_-=c_-a_-$,
$$T_{s_d}^{\underline{j}}=q^{t(a+a_-)/2+t}c_+
J_{s_d}-X_d
{(q^{t(a+a_-)/2+t}-1)X_d-q^{t(a_-+1)/2}+q^{t(a+1)/2}\over
1-X_d^{2}},$$ et, sinon,
$$T_{s_d}^{\underline{j}}=q^{t(a+a_-)/2+t}c_+X_dJ_s-X_d
{(q^{t(a+a_-)/2+t}-1)X_d-q^{t(a_-+1)/2}+q^{t(a+1)/2}\over
1-X_d^{2}}.$$

Notons finalement $M^{\psi }$ le sous-groupe de $M$ engendr\'e par
les \'el\'ements $m$ dont la projection sur chaque facteur $\GL
_{k_{i'}}(F)$ est une puissance $t_{i'}$\`eme, et $B_{\psi }$ la
sous-alg\`ebre de $B$ engendr\'ee par les \'el\'ements $b_m$ avec
$m\in M^{\psi }$. (Dans les notations de \cite{H}, $B_{\psi }$
s'identifie \`a $B_{\so }$ avec $\o $ \'egal \`a l'orbite
inertielle de $\sigma $.)

\null{\bf Proposition:} \it (cf. \cite{H, 5.4, 7.4, 7.6})
L'alg\`ebre $B_{\psi }$, les op\'erateurs
$T_{s_j}^{\underline{j}}$, ainsi que l'op\'erateur
$J_{r_d}^{\underline{j}}$ dans le cas (II), sont contenus dans
$\End _G(i_P^G($ $ind_{M^1}^ME^{\underline{j}}))$.

On a les relations suivantes:

(i) pour $j=1,\dots ,d-1$ et, dans les cas (II) et (IIb), aussi
pour $j=d$,
$$(T_{s_j}^{\underline{j}}+1)(T_{s_j}^{\underline{j}}-q^t)=0.$$ Dans
le cas (III),
$$(T_{s_d}^{\underline{j}}+1)(T_{s_d}^{\underline{j}}-q^{t{{a+a_-}\over
2}+t})=0.$$

(ii) pour $j=1,\dots ,d-2$ et, dans les cas (II) et (IIb) aussi
pour $j=d-1$, on a
$$T_{s_j}^{\underline{j}}T_{s_{j+1}}^{\underline{j}}T_{s_j}^{\underline{j}}
=T_{s_{j+1}}^{\underline{j}}T_{s_j}^{\underline{j}}T_{s_{j+1}}^{\underline{j}}.$$
Dans le cas (III), on trouve
$$T_{s_{d-1}}^{\underline{j}}T_{s_d}^{\underline{j}}T_{s_{d-1}}^{\underline{j}}
T_{s_d}^{\underline{j}}=T_{s_d}^{\underline{j}}T_{s_{d-1}}^{\underline{j}}
T_{s_d}^{\underline{j}}T_{s_{d-1}}^{\underline{j}}.$$ En
particulier, lorsque $w=s_{j_1}\cdots s_{j_l}$ est un \'el\'ement
de $W_{\psi ,i}$ avec $l$ minimal, l'op\'erateur
$T_{s_{j_1}}^{\underline{j}}\cdots T_{s_{j_l}}^{\underline{j}}$ ne
d\'epend que de $w$ et non pas de la d\'ecomposition de $w$ en
sym\'etrie simple choisie.

(iii) Soit $m$ la puissance $t$\`eme d'un \'el\'ement de $M$.
Alors, pour $j=1,\dots ,d-1$ et, dans les cas (II) et (IIb),
\'egalement pour $j=d$,
$$b_mT_{s_j}^{\underline{j}}-T_{s_j}^{\underline{j}}b_{s_j(m)}=
(q^t-1){b_m-b_{s_j(m)}\over 1-X_j^{-1}}.$$ Dans le cas (III), on
trouve
$$b_mT_{s_d}^{\underline{j}}-T_{s_d}^{\underline{j}}b_{s_d(m)}=(q^{t{a+a_-\over
2}+t}-1+X_d^{-1}(q^{t{a+1\over 2}}-q^{t{a_-+1\over
2}})){b_m-b_{s_d(m)}\over 1-X_d^{-2}}$$

(iv) Dans le cas (II) finalement, on a de plus, pour $j=1,\dots
,d$, $T_{s_j}^{\underline{j}}J_{r_d}^{\underline{j}}=
J_{r_d}^{\underline{j}}T_{r_d^{-1}s_jr_d}^{\underline{j}}$,
$(J_{r_d}^{\underline{j}})^2$ est un op\'erateur scalaire et,
lorsque $b\in B_{\psi }$, $J_{r_d}^{\underline{j}}b=\
^{r_d}bJ_{r_d}^{\underline{j}}$. \rm

\null{\bf 3.5} La proposition ci-dessus nous conduit \`a
d\'efinir, pour $w=s_1\cdots s_l$ dans $W_{\psi ,i}$,
$T_w=T_{s_1}^{\underline{j}}\cdots T_{s_l}^{\underline{j}},$
lorsque $s_1\cdots s_l$ est une d\'ecomposition minimale de $w$.

\null {\bf 4.} On part toujours de la situation donn\'ee dans {\bf
2.1}, mais on ne fixe plus $i$. L'objet de section est de
d\'ecrire $\End _G (i_P^G(ind_{M^1}^M E^{\underline{j}}))$ pour un
$\underline{j}$ dans $\J$ fix\'e. Pour cela, notons $W_{\psi
}^{\circ }$ le sous-groupe de $W_{\psi }$ engendr\'e par les
\'el\'ements $s_{i,j}$, posons $r_i=r_{i,d_i}$, lorsque $\rho _i$
v\'erifie les hypoth\`eses du cas (II), et notons $R_{\psi }$ le
sous-groupe de $W_{\psi }$ engendr\'e par ces $r_i$. Le groupe
$W_{\psi }$ est alors le produit semi-direct de $R_{\psi }$ avec
le sous-groupe normal $W_{\psi }^{\circ }$ \cite{H, 1.12}. On
\'ecrira $W_{\psi ,i}^{\circ }$ pour l'intersection $W_{\psi
}^{\circ }\cap W_{\psi ,i}$ et $R_{\psi ,i}$ pour $\{1,r_i\}$.

\null{\bf 4.1 Proposition:} \it (i) Si $w\in W_{\psi ,i}^{\circ }$
et $w'\in W_{\psi ,i'}^{\circ }$ avec $i\ne i'$, alors les
op\'erateurs $T_w^{\underline{j}}$ et $T_{w'}^{\underline{j}}$
commutent.

(ii) Si $r\in R_{\psi ,i}$ et $w\in W_{\psi ,i'}^{\circ }$ avec
$i\ne i'$, alors les op\'erateurs $J_r^{\underline{j}}$ et
$T_w^{\underline{j}}$ commutent.

(iii) Pour tout $r,r'$ dans $R_{\psi }$, les op\'erateurs
$J_r^{\underline{j}}$ et $J_{r'}^{\underline{j}}$ commutent. \rm

\null La proposition nous le permet de d\'efinir
$T_w^{\underline{j}}$ pour tout $w\in W_{\psi }^{\circ }$ et
$J_r^{\underline{j}}$ pour $r\in R_{\psi }$, $r\ne 1$, en posant
$T_w^{\underline{j}}= T_{w_1}^{\underline{j}}\cdots
T_{w_h}^{\underline{j}}$, lorsque $w=w_1\cdots w_h$ avec $w_i\in
W_{\psi _i}$, et $J_r^{\underline{j}}
=J_{r_1}^{\underline{j}}\cdots J_{r_h}^{\underline{j}}$, lorsque
$r=r_1\cdots r_h$ avec $r_i\in R_{\psi ,i}$.

\null {\bf 4.2 Th\'eor\`eme:} \it (cf. \cite{H, 7.7}) Soit
$\underline{j}$ dans $\J $. Les op\'erateurs
$J_r^{\underline{j}}T_w^{\underline{j}}$, $w\in W_{\psi }^{\circ
}$ et $r\in R_{\psi }$, forment une base du $B_{\psi }$-module
$\End _G (i_P^G(ind_{M^1}^M E^{\underline{j}}))$.

La sous-alg\`ebre de $\End _G (i_P^G(ind_{M^1}^M
E^{\underline{j}}))$ engendr\'ee par $B_{\psi }$ et les
$T_w^{\underline{j}}$ est une alg\`ebre de Hecke avec
param\`etres. En particulier, $\End _G (i_P^G(ind_{M^1}^M
E^{\underline{j}}))$ est isomorphe au produit semi-direct de cette
alg\`ebre de Hecke avec param\`etres avec un sous-groupe isomorphe
\`a $R_{\psi }$. \rm

\null\it Remarque: \rm On a $R_{\psi }\ne 1$, si et seulement si
au moins un des $\rho _i$ n'appara\^\i t pas dans $Jord(\psi _H)$
et que $\rho _i$ et $\widehat{G}$ sont ou tous les deux
orthogonaux ou tous les deux symplectiques, ainsi que, dans le cas
$G$ d\'eploy\'e et orthogonal pair, ou bien $k_i$ pair ou bien
$H\ne 1$ et $\tau $ invariant par l'automorphisme ext\'erieur de
$H$.

\null {\bf 5.} On va maintenant proc\'eder \`a la description de
$\End _G(i_P^GE_B)$, lorsque $(\sigma ,E)$ est une
repr\'esentation irr\'eductible cuspidale d'un sous-groupe de Levi
$M$ de $G$ de param\`etre de Langlands $(\psi ,\epsilon )$ (voir
{\bf 1.2} pour la d\'efinition de ce param\`etre, lorsque $G$ est
symplectique ou orthogonal impair, et la section {\bf 1.4},
lorsque $G$ est orthogonal pair). A \'equivalence pr\`es, on peut
supposer que $M$ s'identifie \`a un produit
$$\GL_{k_1}(F)\times\cdots\GL_{k_1}(F)\times \GL_{k_2} (F) \times
\cdots \times\GL_{k_h}\times\cdots\times\GL_{k_h}(F) \times H,$$
chaque facteur $\GL _{k_i}$ \'etant r\'ep\'et\'e $d_i$ fois (cf.
{\bf 2.1}) et que $\psi $ soit de la forme $$\ti{\rho
}_{1,1}\otimes\cdots\ti{\rho }_{1,d_1}\otimes\ti{\rho
}_{2,1}\otimes\cdots\otimes\ti{\rho
}_{2,d_2}\otimes\cdots\otimes\ti{\rho }_{h,d_h}\otimes\psi _H,$$
o\`u les $\ti{\rho }_{i,j}:W_F\rightarrow\GL _{k_i}(\Bbb C)$ sont
des repr\'esentations de $W_F$ qui sont, pour $i$ fix\'e, la
tordue d'une m\^eme repr\'esentation irr\'eductible cuspidale
unitaire $\rho _i$ par un caract\`ere non ramifi\'e. On peut (et
on va) choisir $\rho _i$ autoduale si la contragr\'ediente de
$\ti{\rho }_{i,j}$ est isomorphe au produit de $\ti{\rho }_{i,j}$
par un caract\`ere non ramifi\'e, et on supposera que $\rho _i$ ne
soit pas la tordue d'une $\rho _j$, $j\ne i$, par un caract\`ere
non ramifi\'e.

Rappelons que l'on a d\'esign\'e par $\rho _{i-}$ la
repr\'esentation de $W_F$, d\'etermin\'ee \`a isomorphisme pr\`es,
qui est le produit de $\rho _i$ par un caract\`ere non ramifi\'e
et qui est autoduale et non isomorphe \`a $\rho _i$. \'Ecrivons
$a_i$ (resp. $a_{i-}$) pour l'entier $a_{\rho _i,\psi }$ (resp.
$a_{\rho _{i-},\psi }$) d\'efini par C. Moeglin (cf. {\bf 1.2 -
1.4}) et $t_i$ pour l'ordre du groupe des caract\`eres non
ramifi\'es de $\GL_{k_i}(F)$, stabilisant la classe d'isomorphie
de $\rho _i$. On peut (et on va) supposer $a_i \geq a_{i-}$.

Rappelons finalement que $i_P^GE_B=\bigoplus
_{\underline{j}\in\Cal J}i_P^G(\ind_{M^1}^ME^{\underline{j}})$,
l'ensemble $\Cal J$ et les espaces $E^{\underline{j}}$ ayant
\'et\'e d\'efinis dans {\bf 2.4}.

\null{\bf 5.1} On a d\'efini dans {\bf 2.2}, {\bf 3.1} et {\bf 4.}
le groupe de Weyl $W_{\psi }=W_{\psi }^{\circ }\rtimes R_{\psi }$
qui est engendr\'e par les sym\'etries simples $s_{i,j}$ de
$W_{\psi }^{\circ }$ et $R_{\psi }$. On avait d\'esign\'e dans
{\bf 3.4} par $M^{\psi }$ le sous-groupe de $M$ engendr\'e par les
\'el\'ements $m$ dont la projection sur chaque facteur $\GL
_{k_{i'}}(F)$ est une puissance $t_{i'}$\`eme, et $B_{\psi }$ la
sous-alg\`ebre de $B$ engendr\'ee par les \'el\'ements $b_m$ avec
$m\in M^{\psi }$. Dans {\bf 4.2} on a d\'efini, pour tout
$j\in\Cal J$, $w\in W_{\psi }^{\circ }$ et $r\in R_{\psi }$, des
op\'erateurs $T_w^{\underline{j}}$ et $J_r^{\underline{j}}$ dans
$\End _G(i_P^G (\ind_{M^1}^M E^{\underline{j}}))$. Posons, pour
$w\in W_{\psi }^{\circ }$ et $r\in R_{\psi }$,
$$T_w=\bigoplus _{\underline{j}} T_w^{\underline{j}}\qquad\hbox{\rm et}
\qquad J_r=\bigoplus _{\underline{j}}J_r^{\underline{j}}.$$
\'Ecrivons $X_{i,j}$ pour l'\'el\'ement $X_j$, d\'efini dans {\bf
3.4} relatif \`a $\rho _i$ et $\psi _H$. Pour $\underline{j}\in\J
$, notons $h^{\underline{j}}$ l'\'el\'ement de $M$ qui est le
produit des $(h_{k_i}^{(i,j)})^{j_{i,j}}$ (cf. {\bf 3.2}), $1\leq
i\leq h$, $1\leq j\leq d_i$. Si $\chi\in \Stab(\o )$,
l'application $\phi _{\chi }:i_P^GE_B\rightarrow i_P^GE_B$ qui
agit sur chaque sous-espace $i_P^G(ind_{M^1}^ME^{\underline{j}})$
par le scalaire $\chi (h^{\underline{j}})$ est un automorphisme
$G$-\'equivariant de $i_P^GE_B$ \cite{H, 2.6}.

\null{\bf 5.2 Th\'eor\`eme:} \it L'ensemble des $\phi _{\chi
}J_rT_w$ avec $\chi\in\Stab (\o )$, $r\in R_{\psi }$ et $w\in
W^{\circ }_{\psi }$, forme une base du $B$-module $\End
_G(i_P^GE_B)$. La sous-$B^{\psi }$-alg\`ebre engendr\'ee par
$B^{\psi }$ et les op\'erateurs $T_w$ est une alg\`ebre de Hecke
avec param\`etres. En particulier, pour tous $i=1,\dots ,h$ on
trouve:

(i) pour $j=1,\dots ,d_i-1$ et, si $\rho _i$ v\'erifie
l'hypoth\`ese (II) ou (IIb) de {\bf 3.1}, \'egalement pour
$j=d_i$, on a
$$(T_{s_{i,j}}+1)(T_{s_{i,j}}-q^{t_i})=0.$$ Si $\rho _i$ v\'erifie
l'hypoth\`ese (III), on a
$$(T_{s_{i,d_i}}+1)(T_{s_{i,d_i}}-q^{t_i{a_i+a_{i-}\over
2}+t_i})=0.$$

(ii) pour $j=1,\dots ,d_i-2$ et, si $\rho _i$ v\'erifie
l'hypoth\`ese (II) ou (IIb), \'egalement pour $j=d_i-1$, on a
$$T_{s_{i,j}}T_{s_{i,j+1}}T_{s_{i,j}}
=T_{s_{i,j+1}}T_{s_{i,j}}T_{s_{i,j+1}}.$$ Si $\rho _i$ v\'erifie
l'hypoth\`ese (III), on a
$$T_{s_{i,d_i-1}}T_{s_{i,d_i}}T_{s_{i,d_i-1}}T_{s_{i,d_i}}=T_{s_{i,d_i}}T_{s_{i,d_i-1}}
T_{s_{i,d_i}}T_{s_{i,d_i-1}}.$$ En particulier, lorsque
$w=s_{i_1,j_1}\cdots s_{i_l,j_l}$ est dans $W_{\psi}^{\circ }$
avec $l$ minimal, alors l'op\'era-teur $T_{s_{i_1,j_1}}\cdots
T_{s_{i_lj_l}}$ ne d\'epend que de $w$ et non pas de la
d\'ecomposition choisie.

(iii) Soit $b$ dans $B$. Si $j=1,\dots ,d-1$ et, si $\rho _i$
v\'erifie l'hypoth\`ese (II) ou (IIb), \'egalement pour $j=d_i$,
on a
$$bT_{s_{i,j}}-T_{s_{i,j}}\ ^{s_{i,j}}b=(q^{t_i}-1){b-\ ^{s_{i,j}}b\over
1-X_{i,j}^{-1}}.$$ Si $\rho _i$ v\'erifie l'hypoth\`ese (III), on
trouve
$$bT_{s_{i,d_i}}-T_{s_{i,d_i}}\ ^{s_{i,d_i}}b=(q^{t_i{a_i+a_{i-}\over 2}+t_i}-1+X_{i,d_i}^{-1}(q^{t_i
{a_i+1\over 2}}-q^{t_i{a_{i-}+1\over 2}})){b-\ ^{s_{i,d_i}}b\over
1-X_{i,d_i}^{-2}}$$

Par ailleurs,

(iv) Pour tout $r\in R_{\psi }$, $b\in B$ et $s\in W_{\psi }$, on
a $T_sJ_r=J_rT_{r^{-1}sr}$, $J_r^2$ est un op\'erateur scalaire et
$J_rb=\ ^rbJ_r$.

(v) Les automorphismes $\phi _{\chi }$ commutent avec les $T_w$,
$J_r$ et les \'el\'ements de $B_{\psi }$. En g\'en\'eral, on a
pour $m\in M$, $\phi _{\chi }b_m=\chi (m)b_m\phi _{\chi }$.\rm

\null\it Remarque: \rm Notons $Rat(M)$ l'ensemble des caract\`eres
rationnels de $M$. Dans le langage de Lusztig \cite{L},
l'alg\`ebre de Hecke du th\'eor\`eme est associ\'ee \`a la
donn\'ee radicielle basique $(\Lambda ,\Sigma ,\Lambda ^{\vee
},\Sigma ^{\vee },\Delta )$ avec $\Lambda$ \'egal \`a l'image de
$M^{\psi }$ dans $M/M^1$, $\Lambda ^{\vee }$ \'egal au
sous-ensemble de $Rat(M)\otimes\Bbb R$ qui est en dualit\'e
parfaite avec $\Lambda $ par $(m,\chi\otimes x)\mapsto
-x\log_q(\chi (m))$, $\Sigma =\{^wX_{i,j}\vert w\in W_{\psi
}^{\circ },i=1,\dots ,h\ \hbox{\rm et}\ j=1,\dots ,d_i'\}$,
$\Sigma ^{\vee }$ \'egal au syst\`eme de racines dual dans
$\Lambda ^{\vee }$ et $\Delta =\{X_{i,j}\vert i=1,\dots ,h\
\hbox{\rm et}\ j=1,\dots ,d_i'\}$. (Rappelons que $d_i'$ a \'et\'e
d\'efini \`a la fin de {\bf 3.1}. Le fait que ceci forme bien un
syst\`eme de racines est prouv\'e dans \cite{H, 6}.)

Les param\`etres de l'alg\`ebre de Hecke dans le langage de
Lusztig (il propose plusieurs terminologies) se lisent directement
sur les relations donn\'ees dans le th\'eor\`eme {\bf 5.2} selon
les diff\'erents cas. Remarquons que les param\`etres
$t_i{a_i+a_{i-}\over 2}$ sont toujours des entiers suite \`a la
proposition {\bf 1.3}.

On ne peut pas prendre pour $\Lambda $ le r\'eseau $M/M^1$ tout
entier, puisque $\Sigma ^{\vee }$ n'est pas inclus dans le
r\'eseau dual (\`a moins que tous les $t_i$ soient nuls). Donc, la
$B$-sous-alg\`ebre engendr\'ee par les op\'erateurs $T_w$ n'est
pas une alg\`ebre de Hecke avec param\`etres au sens strict de la
terminologie de Lusztig.

\null {\bf 5.3 Lemme:} \it Pour tout $\underline{j}$, la
projection $p^{\underline{j}}$ de $i_P^GE_B$ sur
$i_P^GE^{\underline{j}}$ est une combinaison $\Bbb C$-lin\'eaire
des automorphismes $\phi _{\chi }$, $\chi\in Stab(\o )$.

\null Preuve: \rm Sans perte de g\'en\'eralit\'e, on peut supposer
$\underline{j}=(1,1,\dots ,1)=:\underline{1}.$ Pour
$\underline{j}$ diff\'erent de $\underline{1}$, on peut choisir un
caract\`ere $\chi _{\underline{j}}$ dans $\Stab(\o )$ tel que
$\chi_{\underline{j}} (h^{\underline{j}})\ne 1$. L'endomorphisme
$\phi_{\chi }-\chi_{\underline{j}} (h^{\underline{j}})$ est alors
trivial sur $i_P^GE^{\underline{j}}$ et non trivial sur $i_P^G
E^{\underline{1}}$. Le compos\'e de tous ces endomorphismes pour
$\underline{j}\ne\underline{1}$ (dans n'importe quel ordre) est
alors un endomorphisme de $i_P^GE_B$ qui est trivial sur chaque
sous-espace $i_P^GE^{\underline{j}}$,
$\underline{j}\ne\underline{1}$, et qui agit sur
$i_P^GE^{\underline{1}}$ par un scalaire non nul \'egal au produit
des $1-\chi_{\underline{j}}(h^{\underline{j}}).$

\hfill{\fin 2}

\null {\bf 5.4} \it Preuve du th\'eor\`eme {\bf 5.2}: \rm Il est
clair que $B$, les automorphismes $\phi _{\chi }$ et les $T_w$
ainsi que les $J_r$ sont dans $\End _G(i_P^GE_B)$. Il reste \`a
voir que la sous-alg\`ebre $\A $ engendr\'ee par ces op\'erateurs
est \'egale \`a $\End _G(i_P^GE_B)$. Soit $\Phi\in \End
_G(i_P^GE_B)$. Comme $\Phi =\bigoplus
_{\underline{j}}p^{\underline{j}}\circ\Phi $, on est ramen\'e au
cas o\`u l'image de $\Phi$ est contenue dans un des espaces
$i_P^GE^{\underline{j}}$.

Fixons $\underline{j}$. Pour tout $\underline{i}$ dans $\J$, il
existe un \'el\'ement $m_{\underline{i}}$ de $M$ tel que
$b_{m_{\underline{i}}}i_P^GE^{\underline{j}}=i_P^GE^{\underline{i}}.$
La restriction \`a gauche de $b_{m_{\underline{i}}}\Phi $ \`a
l'espace $i_P^GE^{\underline{i}}$ d\'efinit un endomorphisme de
$i_P^GE^{\underline{i}}$. Par le th\`eor\`eme {\bf 4.2}, c'est
donc une combinaison $B_{\psi }$-lin\'eaire des
$J_r^{\underline{i}}T_w^{\underline{i}}$, $r\in R_{\psi }$ et
$w\in W_{\psi }^{\circ }$. Notons $b_{r,w}^{\underline{i}}$ les
coefficients dans cette combinaison lin\'eaire. On trouve alors
$$\Phi =\bigoplus_{\underline{i}}b_{m_{\underline{i}}}^{-1}\sum
_{r,w}b_{r,w}^{\underline{i}}J_r^{\underline{i}}T_w^{\underline{i}}.$$
Or, comme $T_w^{\underline{i}}=p^{\underline{i}}T_w$ et
$J_r^{\underline{i}} =p^{\underline{i}}J_r$, cet op\'erateur
appartient suite au lemme {\bf 5.3} \`a $\A $. Il en est donc de
m\^eme de $\Phi $.\hfill{\fin 2}

\enddocument
\bye